\documentclass[11pt,twoside,]{amsart}
\usepackage{amsmath, amsthm, amscd, amsfonts, amssymb, color}
\usepackage[dvips]{graphicx}
\usepackage{tikz}
\usepackage[bookmarksnumbered, plainpages]{hyperref}
\newtheorem{thm}{Theorem}[section]

\newtheorem{lem}[thm]{Lemma}
\newtheorem{prop}[thm]{Proposition}

\numberwithin{equation}{section}
\def\pn{\par\noindent}

\begin{document}

\title{Connected graphs cospectral with a Friendship graph}
\author{Alireza Abdollahi$^*$ and  Shahrooz Janbaz}

\thanks{{\scriptsize
\hskip -0.4 true cm MSC(2010): Primary: 05C50; Secondary: 05C31.
$*$Corresponding author}}
\maketitle

\begin{abstract}
Let $n$ be any positive integer, the friendship graph $F_n$ consists of $n$ edge-disjoint triangles that all of them meeting in one vertex. A graph $G$ is called cospectral with a graph $H$ if their adjacency matrices have the same eigenvalues. Recently in \href{http://arxiv.org/pdf/1310.6529v1.pdf}{http://arxiv.org/pdf/1310.6529v1.pdf} it is proved that if $G$ is any graph cospectral with $F_n$ ($n\neq 16$), then $G\cong F_n$. In this note, we give a proof of a special case of the latter: Any connected graph cospectral with $F_n$ is isomorphic to $F_n$.
Our proof is independent of ones given in \href{http://arxiv.org/pdf/1310.6529v1.pdf}{http://arxiv.org/pdf/1310.6529v1.pdf} and the proofs are based on our recent results given in [{\em Trans.  Comb.}, {\bf 2} no. 4 (2013) 37-52.] using an upper bound for the largest eigenvalue of a connected graph given in
[{\em J. Combinatorial Theory Ser. B} {\bf 81} (2001) 177-183.].
\end{abstract}

\vskip 0.2 true cm

\section{\bf Introduction}
\vskip 0.4 true cm

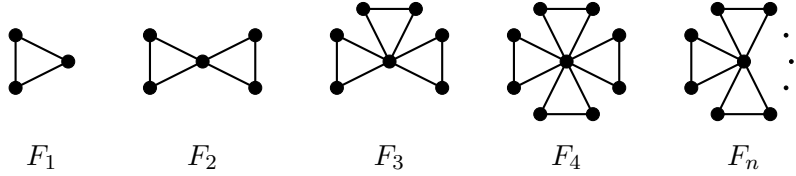
\begin{figure}[htb]
\centering
\begin{tikzpicture}[scale=0.7]
\filldraw [black]
(-1,-0.5) circle (3.5 pt)
(-1,0.5) circle (3.5 pt)
(0,0) circle (3.5 pt);
\node [label=below:$F_1$] (F_1) at (-0.5,-1.2) {};
\draw[thick] (-1,-0.5) -- (-1,0.5);
\draw[thick] (-1,0.5) -- (0,0);
\draw[thick] (-1,-0.5) -- (0,0);
\end{tikzpicture} \qquad
\begin{tikzpicture}[scale=0.7]
\filldraw [black]
(-1,-0.5) circle (3.5 pt)
(-1,0.5) circle (3.5 pt)
(0,0) circle (3.5 pt)
(1,0.5) circle (3.5 pt)
(1,-0.5) circle (3.5 pt);
\node [label=below:$F_2$] (F_2) at (0,-1.2) {};
\draw[thick] (-1,-0.5) -- (-1,0.5);
\draw[thick] (-1,0.5) -- (0,0);
\draw[thick] (-1,-0.5) -- (0,0);
\draw[thick] (0,0) -- (1,0.5);
\draw[thick] (0,0) -- (1,-0.5);
\draw[thick] (1,0.5) -- (1,-0.5);
\end{tikzpicture} \qquad
\begin{tikzpicture}[scale=0.7]
\filldraw [black]
(-1,-0.5) circle (3.5 pt)
(-1,0.5) circle (3.5 pt)
(0,0) circle (3.5 pt)
(1,0.5) circle (3.5 pt)
(1,-0.5) circle (3.5 pt)
(-0.5,1) circle (3.5 pt)
(0.5,1) circle (3.5 pt);
\node [label=below:$F_3$] (F_3) at (0,-1.2) {};
\draw[thick] (-1,-0.5) -- (-1,0.5);
\draw[thick] (-1,0.5) -- (0,0);
\draw[thick] (-1,-0.5) -- (0,0);
\draw[thick] (0,0) -- (1,0.5);
\draw[thick] (0,0) -- (1,-0.5);
\draw[thick] (1,0.5) -- (1,-0.5);
\draw[thick] (-0.5,1) -- (0.5,1);
\draw[thick] (-0.5,1) -- (0,0);
\draw[thick] (0.5,1) -- (0,0);
\end{tikzpicture}\qquad
\begin{tikzpicture}[scale=0.7]
\filldraw [black]
(-1,-0.5) circle (3.5 pt)
(-1,0.5) circle (3.5 pt)
(0,0) circle (3.5 pt)
(1,0.5) circle (3.5 pt)
(1,-0.5) circle (3.5 pt)
(-0.5,1) circle (3.5 pt)
(0.5,1) circle (3.5 pt)
(-0.5,-1) circle (3.5 pt)
(0.5,-1) circle (3.5 pt);
\node [label=below:$F_4$] (F_4) at (0,-1.2) {};
\draw[thick] (-1,-0.5) -- (-1,0.5);
\draw[thick] (-1,0.5) -- (0,0);
\draw[thick] (-1,-0.5) -- (0,0);
\draw[thick] (0,0) -- (1,0.5);
\draw[thick] (0,0) -- (1,-0.5);
\draw[thick] (1,0.5) -- (1,-0.5);
\draw[thick] (-0.5,1) -- (0.5,1);
\draw[thick] (-0.5,1) -- (0,0);
\draw[thick] (0.5,1) -- (0,0);
\draw[thick] (-0.5,-1) -- (0.5,-1);
\draw[thick] (-0.5,-1) -- (0,0);
\draw[thick] (0.5,-1) -- (0,0);
\end{tikzpicture}\qquad
\begin{tikzpicture}[scale=0.7]
\filldraw [black]
(-1,-0.5) circle (3.5 pt)
(-1,0.5) circle (3.5 pt)
(0,0) circle (3.5 pt)
(-0.5,1) circle (3.5 pt)
(0.5,1) circle (3.5 pt)
(0.8,0.5) circle (1 pt)
(0.9,0) circle (1 pt)
(0.8,-0.5) circle (1 pt)
(-0.5,-1) circle (3.5 pt)
(0.5,-1) circle (3.5 pt);
\node [label=below:$F_n$] (F_n) at (0,-1.2) {};
\draw[thick] (-1,-0.5) -- (-1,0.5);
\draw[thick] (-1,0.5) -- (0,0);
\draw[thick] (-1,-0.5) -- (0,0);
\draw[thick] (-0.5,1) -- (0.5,1);
\draw[thick] (-0.5,1) -- (0,0);
\draw[thick] (0.5,1) -- (0,0);
\draw[thick] (-0.5,-1) -- (0.5,-1);
\draw[thick] (-0.5,-1) -- (0,0);
\draw[thick] (0.5,-1) -- (0,0);
\end{tikzpicture}
\caption{Friendship graphs $F_1$, $F_2$, $F_3$, $F_4$ and $F_n$}\label{pic1}
\end{figure}

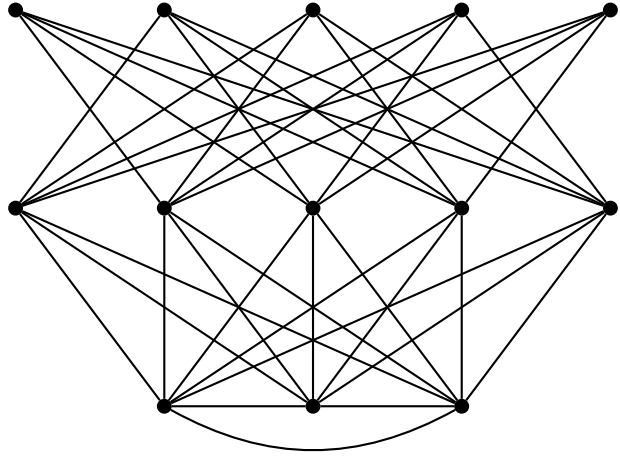
\begin{figure}[htb]
\centering
\usetikzlibrary{shapes.geometric}
\begin{tikzpicture}
[every node/.style={inner sep=0pt}]
\node (1) [circle, minimum size=5pt, fill=black, line width=0.625pt, draw=black] at (112.5pt, -37.5pt)  {};
\node (2) [circle, minimum size=5pt, fill=black, line width=0.625pt, draw=black] at (168.75pt, -37.5pt)  {};
\node (3) [circle, minimum size=5pt, fill=black, line width=0.625pt, draw=black] at (225.0pt, -37.5pt)  {};
\node (4) [circle, minimum size=5pt, fill=black, line width=0.625pt, draw=black] at (281.25pt, -37.5pt)  {};
\node (5) [circle, minimum size=5pt, fill=black, line width=0.625pt, draw=black] at (337.5pt, -37.5pt)  {};
\node (6) [circle, minimum size=5pt, fill=black, line width=0.625pt, draw=black] at (112.5pt, -112.5pt)  {};
\node (7) [circle, minimum size=5pt, fill=black, line width=0.625pt, draw=black] at (168.75pt, -112.5pt)  {};
\node (8) [circle, minimum size=5pt, fill=black, line width=0.625pt, draw=black] at (225.0pt, -112.5pt)  {};
\node (9) [circle, minimum size=5pt, fill=black, line width=0.625pt, draw=black] at (281.25pt, -112.5pt)  {};
\node (10) [circle, minimum size=5pt, fill=black, line width=0.625pt, draw=black] at (337.5pt, -112.5pt)  {};
\node (11) [circle, minimum size=5pt, fill=black, line width=0.625pt, draw=black] at (168.75pt, -187.5pt)  {};
\node (12) [circle, minimum size=5pt, fill=black, line width=0.625pt, draw=black] at (225.0pt, -187.5pt)  {};
\node (13) [circle, minimum size=5pt, fill=black, line width=0.625pt, draw=black] at (281.25pt, -187.5pt)  {};
\draw [line width=0.8, color=black] (1) to  (7);
\draw [line width=0.8, color=black] (1) to  (8);
\draw [line width=0.8, color=black] (1) to  (9);
\draw [line width=0.8, color=black] (1) to  (10);
\draw [line width=0.8, color=black] (2) to  (6);
\draw [line width=0.8, color=black] (2) to  (8);
\draw [line width=0.8, color=black] (2) to  (9);
\draw [line width=0.8, color=black] (2) to  (10);
\draw [line width=0.8, color=black] (3) to  (6);
\draw [line width=0.8, color=black] (3) to  (7);
\draw [line width=0.8, color=black] (3) to  (9);
\draw [line width=0.8, color=black] (3) to  (10);
\draw [line width=0.8, color=black] (4) to  (6);
\draw [line width=0.8, color=black] (4) to  (7);
\draw [line width=0.8, color=black] (4) to  (8);
\draw [line width=0.8, color=black] (4) to  (10);
\draw [line width=0.8, color=black] (5) to  (6);
\draw [line width=0.8, color=black] (5) to  (7);
\draw [line width=0.8, color=black] (5) to  (8);
\draw [line width=0.8, color=black] (5) to  (9);
\draw [line width=0.8, color=black] (11) to  (12);
\draw [line width=0.8, color=black] (12) to  (13);
\draw [line width=0.8, color=black] (11) to  [in=209, out=331] (13);
\draw [line width=0.8, color=black] (11) to  (6);
\draw [line width=0.8, color=black] (11) to  (7);
\draw [line width=0.8, color=black] (11) to  (8);
\draw [line width=0.8, color=black] (11) to  (9);
\draw [line width=0.8, color=black] (11) to  (10);
\draw [line width=0.8, color=black] (12) to  (6);
\draw [line width=0.8, color=black] (12) to  (7);
\draw [line width=0.8, color=black] (12) to  (8);
\draw [line width=0.8, color=black] (12) to  (9);
\draw [line width=0.8, color=black] (12) to  (10);
\draw [line width=0.8, color=black] (13) to  (10);
\draw [line width=0.8, color=black] (13) to  (9);
\draw [line width=0.8, color=black] (13) to  (8);
\draw [line width=0.8, color=black] (13) to  (7);
\draw [line width=0.8, color=black] (13) to  (6);
\end{tikzpicture}
\caption{The disjoint union of the above graph with $10K_2$ is cospectral with $F_{16}$}\label{pic2}
\end{figure}

The friendship graph $F_n$ is the graph consisting of $n$ edge-disjoint cycles of length $3$, all meeting in a common vertex (see Figure \ref{pic1}). The adjacency spectrum $Spec(G)$ of a graph $G$ is the multiset of eigenvalues of its adjacency matrix. A graph $G$ is called determined by  its adjacency spectrum (for short $DS$), if $Spec(G)=Spec(H)$ for some graph $H$, then $G\cong H$.  In  \cite{WangB,Wang}, it is conjectured that the friendship graph is $DS$.  In  \cite{Das}, it is claimed that the conjecture is valid, but  it is noted in \cite{Ab-Ja} that the proof has a flaw. In \cite{Ab-Ja} some partial results on graphs cospectral with $F_n$ has been obtained. Finally in \cite{Ham} a more general result about
graphs with all but two eigenvalues equal to $\pm 1$ is obtained and as a corollary it is shown that $F_n$ is DS if $n \neq 16$ and $F_{16}$ is cospectral with the disjoin union of $10$ complete graph $K_2$ with the graph given in  Figure \ref{pic2}. The spectrum of the graph of Figure \ref{pic2} is $\left\lbrace \left(1\pm \sqrt{129}\right)/2,1^5,-1^6 \right\rbrace$. 
Our main result is the following:
\begin{thm}\label{main}
Any connected graph cospectral with $F_n$ is isomorphic to $F_n$.
\end{thm}
Our proof is based on our previous recent results in \cite{Ab-Ja} by using an upper bound given in \cite{Hong} for the spectral radius of a connected graph in terms of the minimum degree of the graph and its number of vertices and edges.  Our proof is  independent of that of  given in \cite{Ham}.  

\noindent

\noindent


\section{\bf Proof of Theorem \ref{main}}
\vskip 0.4 true cm
In the following, we give results that we will need in our proof of Theorem \ref{main}. 
\begin{prop}\cite[Proposition 2.3]{Ab-Ja}\label{Ab-Ja-1}
Let $F_n$ denote the friendship graph with $2n+1$ vertices. Then
$$Spec(F_n)=\left\lbrace\left(\frac{1}{2}-\frac{1}{2}\sqrt{1+8n}\right)^1,{\left(-1\right)}^n,{\left(1\right)}^{n-1},{\left(\frac{1}{2}+\frac{1}{2}\sqrt{1+8n}\right)}^1 \right\rbrace.$$
\end{prop}

The maximum eigenvalue of a graph $G$ is called spectral radius and it is denoted by $\varrho(G)$.
\begin{lem}\cite[Lemma 3.2]{Ab-Ja}\label{2}
Let $G$ be a connected graph that is cospectral with $F_n$ and $\delta(G)$ be the minimum degree of $G$. Then, $\delta(G)=2$ and $G$ has at least $1+\varrho(F_n)$ vertices with this minimum degree.
\end{lem}


\begin{thm}\cite[Theorem 3.11]{Ab-Ja}\label{Ab-Ja-4}
Let $G$ be a graph cospectral with $F_n$ and $G$ has two adjacent vertices of degree $2$. Then $G$ is isomorphic to $F_n$.
\end{thm}

\begin{thm}\cite[Theorem 2.3]{Hong}\label{Hong}
Let $G$ be a simple graph with $n$ vertices and $m$ edges. Let $\delta=\delta(G)$ be the minimum degree of vertices of $G$ and $\varrho (G)$ be the spectral radius of the adjacency matrix  of $G$. Then
$$\varrho (G)\leq \frac{\delta -1}{2}+\sqrt{2m-n\delta+\frac{\left(\delta +1 \right)^2}{4}}.$$
Equality holds if and only if $G$ is either a regular graph or a bidegreed graph in which each vertex is of degree either $\delta$ or $n-1$.
\end{thm}

Now we are ready to prove the main result of this note.\\

\noindent{\bf Proof of Theorem \ref{main}.} Let $G$ be a connected graph and cospectral with the friendship graph $F_n$. 
By Proposition \ref{Ab-Ja-1}, $\varrho(G)=\frac{1}{2}+\frac{1}{2}\sqrt{1+8n}$ and the number vertices and edges of $G$ are the same as $F_n$. Also by Lemma \ref{2}, $\delta=\delta(G)=2$. Therefore, the equality in Theorem \ref{Hong} holds and so  $G$ is either a regular graph or a bidegreed graph in which each vertex is of degree either $2$ or $2n$. If $n=1$, then it is easy to see that $G\cong F_1=K_3$. If $n>1$ then $F_n$ is not regular and so $G$ is not a regular graph, since regularity can be determined by the adjacency spectrum of a graph.  Hence every vertex of $G$ has degree $2$ or $2n$. 
Since $G$ has $2n+1$ vertices and $3n$ edges, it follows that $G$ has only one vertex of degree $2n$ and all other vertices has degree $2$. Thus $G$ has at least two adjacent vertices of degree $2$. Now Theorem \ref{Ab-Ja-4} completes the proof. $\hfill \Box$



\bigskip

{\footnotesize \pn{\bf Alireza Abdollahi}\; \\ {Department of Mathematics, University of Isfahan,} {Isfahan 81746-73441, Iran}\\
and School of Mathematics, Institute for Research in Fundamental Sciences (IPM), P.O.Box: 19395-5746, Tehran, Iran.\\
{\tt Email: a.abdollahi@math.ui.ac.ir}\\

{\footnotesize \pn{\bf Shahrooz Janbaz}\; \\ {Department of Mathematics}, {University of Isfahan,} {Isfahan 81746-73441, Iran}\\
{\tt Email: shahrooz.janbaz@sci.ui.ac.ir}\\

\end{document}